\documentstyle[12pt]{article}

\begin{document}
\newcommand{\ol }{\overline}
\newcommand{\ul }{\underline }
\newcommand{\ra }{\rightarrow }
\newcommand{\lra }{\longrightarrow }
\newcommand{\ga }{\gamma }
\newcommand{\st }{\stackrel }

\title{\Large\bf Higher Schur-multiplicator of a Finite Abelian Group}
\author{by\\M.R.R.Moghaddam and B.Mashayekhy\\ Department of Mathematics \\
Ferdowsi University of Mashhad\\ P.O.Box 1159-91775, Mashhad, Iran\\e-mail:
          moghadam@science2.um.ac.ir\\  mashaf@science2.um.ac.ir }
\date{ }
\maketitle
\begin{abstract}
 In this paper we obtain an explicit formula for the {\it higher
Schur-multiplicator }of an arbitrary finite abelian group with
respect to the
variety of nilpotent groups of class at most $c\geq 1$ .   \\
1991 Mathematics Subject Classification: Primary 20F12,20F18,20K25 \\
Key words and phrases: Abelian Group , Higher Schur-multiplicator
, Nilpotent Variety.
\end{abstract}
\vspace{20mm}
\newpage
\begin{center}{\large\bf 1. Introduction and Preliminaries} \end{center}

 In 1907 I.Schur [7] proved that the Schur-multiplicator of a direct product
of two finite groups is isomorphic to the direct sum of the Schur-multiplicators
of the direct factors and the tensor product of the two groups abelianized.\\
(see also J.Wiegold [8].) If
$$ G={\bf Z}_{n_1}\oplus {\bf Z}_{n_2}\oplus \ldots {\bf Z}_{n_k} $$
is a finite abelian group , where $n_{i+1}|n_i$ , for all $1\leq i\leq k-1$ ,
then using the above fact one obtains the Schur-multiplicator of $G$ as follows
( see [4]):
$$ M(G)\cong {\bf Z}_{n_2}\oplus {\bf Z}^{(2)}_{n_3}\oplus \ldots \oplus
{\bf Z}_{n_k}^{(k-1)}\ \ \ , $$
where ${\bf Z}^{(m)}_{n}$ denotes the direct product of $m$ copies of the cyclic
group ${\bf Z}_n$ .

 Now, in this paper, a similar result will be presented for the {\it higher
Schur-multiplicator} of an arbitrary finite abelian group with respect to the
variety of nilpotent groups of class at most $c\geq 1$ , ${\cal N}_c$ , say.
(see [5] for the notation.)

 Let $G$ be any group with a free presentation
 $$ 1\lra R\lra F\lra G\lra 1 \ \ \ ,$$
where $F$ is a free group. Then, following Heinz Hopf [3],
$$ \frac {R\cap F'}{[R,F]} \ \ ,$$
is isomorphic to the Schur-multiplicator of $G$ , denoted by $M(G)$ . Now , the
{\it higher Schur-multiplicator} of $G$ with respect to the variety of nilpotent
groups of class at most $c\geq 1$ , is defined to be
$$ {\cal N}_cM(G)=\frac {R\cap \ga_{c+1}(F)}{[R,\ _cF]}\ \ \ \ ,$$
where $\ga_{c+1}(F)$ is the $(c+1)$st-term of the lower central series and
$[R,\ _cF]$ denotes $[R,\underbrace {F,F,\ldots ,F}_{c-times}] $ , see [5 or 6]
for further details. \\

 Let $H_i=<x_i|x_i^{r_i}>\cong {\bf Z}_{r_i}$ , $i=1,2,\ldots ,t\ \ ,\ \ r_i\in
{\bf N}$ be cyclic groups of order $r_i$ , $1\leq i\leq t$ , $r_i\geq 0$ and let
$$ 1\lra R_i=<x_i^{r_i}>\st {\pi_i}{\lra }F_i=<x_i>\lra H_i\lra 1 $$
be the free presentation for $H_i$ , where $1\leq i\leq t$ . Also, let
   $$ G=\prod_{i=1}^{t}\!^{\times }H_i={\bf Z}_{r_1}\times {\bf Z}_{r_2}\times
   \ldots \times {\bf Z}_{r_t} $$
be the direct product of cyclic groups ${\bf Z}_{r_i}$'s. Then
$$ 1\lra R\lra F\lra G\lra 1 $$
 is the free presentation for $G$, where
  $$ F=\prod_{i=1}^{t}\!^{*}F_i=<x_1,\ldots ,x_t>\ \  and\ \
R=<x_i^{r_i},\ga_2(F)\cap [F_i]^*|i=1,\ldots ,t>^F\ ,$$
where $\prod_{i=1}^{t}\!^{*}F_i$ is the free product of $F_i$'s , $i=1,2,\ldots
,t$ , and $[F_i]^{*}$ is the normal closure of some commutator subgroups of the
the free product , defined as follows:
 $$[F_i]^{*}=<[F_i,F_j]\ |\ 1\leq i,j\leq t\ ,\ i\neq j\ >^F\ \ \ . $$
Since $F_i$'s are cyclic, we have $\ga_2(F)\subseteq[F_i]^*$ . Hence
$$ R=<x_i^{r_i},\ga_2(F)\ |\ i=1,\ldots ,t>^F\ .$$
 Now, put $S=<x_1^{r_1},\ldots ,x_t^{r_t}>^F$ and for all $m\geq
1$, define $\rho_1(S)=S\ \ ,\ \ \rho_m(S)=[S,\ _{m-1}F]$, inductively. This
yeilds the central series
$$ S=\rho_1(S)\supseteq \rho_2(S)\supseteq \rho_3(S)\supseteq \ldots \supseteq
\rho_m(S)\supseteq \ldots \ \ \ .$$
Thus we have
$$R=S\ga_2(F)\ \ \ {\rm and}\ \ \  [R,\ _mF]=\rho_{m+1}(S)\ga_{m+2}(F)\ \ \ \
.\ (*)   $$

 Let $F=\prod_{i=1}^{t}\!^{*}F_i$ be the free product of $F_1,F_2,\ldots ,F_t$ .
We define a {\it basic commutator subgroup} , $B(F_1,F_2,\ldots ,F_t)_s $ of
weight $s\ (\ s\in {\bf N}\ )$ on $t$ free groups $F_1,F_2,\ldots ,F_t$ ,
as follows:\\
 We first order the subgroups $F_1,F_2,\ldots ,F_t$ by setting
$F_i<F_j$ if $i<j$ . Then $B(F_1,F_2,\ldots ,F_t)_s$ is the subgroup generated
by all the basic commutators of weight $s$ on $t$ letters $x_1,x_2,\ldots ,x_t$
, where $x_\in F_i$ for all $1\leq i\leq t$ . For the definition of basic
commutators see M.Hall [1].

 Note that here we have slightly modified the definition of basic commutator
subgroups from M.R.R.Moghaddam [6].

 Now, let $T(H_1,H_2,\ldots ,H_t)_{s}$ denote the summation of all the tensor
products corresponding to the basic commutator subgroups $B(F_1,F_2,\ldots
,F_t)_{s}$ where
$$ 1\lra R_i\lra F_i\st {\pi_i}{\lra }H_i\lra 1 $$
is the free presentation for $H_i\ ,\ i=1,2,\ldots ,t$ . More precisely, if
$[F_j.F_i,\ldots ]$,\\ with any bracketing, is a basic commutator subgroup of
weight $s$ in $F_i$'s, then the {\it ``corresponding'' } tensor product will be
$$ (H_j\otimes H_i\otimes \ldots )\ \ , $$
bracketed in the same way . (Note that $H_i$'s are abelian groups.)

 Similarly, the element $[x_j,x_i,\ldots ]$ of the commutator subgroup
$B(F_1,F_2,\ldots ,F_t)_s $,\\ with any bracketing , corresponds to
the element of the tensor product
$$( {\pi_jx_j}\otimes  {\pi_ix_i}\otimes \ldots )$$
bracketed in the same way, where $x_k$ is the generator of $F_k$ .

 We keep this notation throughout the rest of the paper , and it will be used
without further reference.

\begin{center}           {\large\bf 2. The Main Results } \\
\end{center}
{\bf Lemma 2.1}

Let $G$ be a finite abelian group , then by the previous notation for all
$c\geq 1$ ,
$$ {\cal N}_cM(G)\cong \frac {\ga_{c+1}(F)}{\rho_{c+1}(S)\ga_{c+2}(F)}\ \ .$$
\\ \ \ \\
{\bf Proof.}

 Let $1\lra R\lra F\lra G\lra 1$ be a free presentation for $G$ , then by the
definition and using the fact that $\ga_{c+1}(F)$ is contained in $S\ga_2(F)=R$
, and $(*)$ ,
$$ {\cal N}_cM(G)\cong \frac {R\cap \ga_{c+1}(F)}{[R,\ _cF]} $$
$$  \ \ \ \ \ \ \ \  \cong \frac {\ga_{c+1}(F)}{[R,\ _cF]}  $$
$$ \ \ \ \ \ \ \ \ \ \ \ \ \ \ \ \ \ \ \ \ \ \ \ \ \ \ \ \cong \frac
{\ga_{c+1}(F)}{\rho_{c+1}(S)\ga_{c+2}(F)} \ \ \ \ \
.\ \ \Box $$

 Now, we are in a position to prove the following important theorem.\\
{\bf Theorem 2.2}

 Consider the above assumption and notation . Then , for all $c\geq 1$ ,
$$ \frac {\ga_{c+1}(F)}{\rho_{c+1}(S)\ga_{c+2}(F)}\cong T(H_1,H_2,\ldots
,H_t)_{c+1}\ \ . $$
{\bf Proof.}

 Since $F_i$'s are infinite cyclic groups , using Hall's Theorem [2] we obtain
$$B(F_1,F_2,\ldots ,F_t)_{c+1}=\ga_{c+1}(F) \ \ \ \ \ \ \
modulo\ \ga_{c+2}(F)\ \ .$$
By the notation of previous
section , there exists an obvious relation between $B(F_1,F_2,\ldots ,F_t)_{c+1}$
and
$T(H_1,H_2,\ldots ,H_t)_{c+1}$ . Also the commutators of weight $c+1$ in
$B(F_1,F_2,\ldots ,F_t)_{c+1}\ \ modulo\ \ga_{c+2}(F)$ are multilinear. Now
one may construct a homomorphism $\mu $ from
$\ga_{c+1}(F)/\rho_{c+1}(S)\ga_{c+2}(F)$
into the abelian group $T(H_1,H_2,\ldots ,H_t)_{c+1}$ given by
$$ \prod \underbrace {[f,g,\ldots
]}_{wt.c+1}\rho_{c+1}(S)\ga_{c+2}(F)\longmapsto
\sum ({\pi_jf},{\pi_ig},\ldots ) $$
$$ with\ any\ bracketing\ \ \ \ \ \ \ \ \ \ \ \ \ {\rm ``}corresponding{\rm "}\ element\ . $$
We know that $\ga_{c+1}(F)/\ga_{c+2}(F)$ is the free abelian group on the basic
commutators of weight $c+1$ on $t$ letters , (by P.Hall [2]) .
Now we define the same mapping, as the above
correspondence, on the set of basic commutators of weight $c+1$ on $t$ letters
into the
abelian group $T(H_i,H_2,\ldots ,H_t)_{c+1}$ . Then by the universal property
of free abelian groups it can be extended to a homomorphism $\phi $, say, from
$\ga_{c+1}(F)/\ga_{c+2}(F)$ into $T(H_1,H_2,\ldots ,H_t)_{c+1} $ .

 To show that $\phi $ induces the homomorphism $\mu $ from
$\ga_{c+1}(F)/\rho_{c+1}(S)\ga_{c+2}(F)$ into $T(H_1,H_2,\ldots ,H_t)_{c+1}$ ,
we must prove that $\rho_{c+1}(S)$ is mapped onto zero under $\phi $ . But this
is obvious, since
$$\underbrace{[x_i^{r_i},f,g,\ldots ]}_{wt.c+1}\equiv
\underbrace{[x_i,f,g,\ldots ]^{r_i}}_{wt.c+1}\ \ \ \ (mod\ \ga_{c+2}(F)) $$
and
$$ [x_i,f,g,\ldots ]^{r_i}\st {\phi }{\longmapsto }r_i({\pi_ix}\otimes
{\pi_jf}\otimes {\pi_kg}\otimes \ldots ) $$
$$\ \ \ \ \ \ \ \ \ \ \ \ \ \ \ \ \ \ \ \ \ \ \ \ \  =(r_i{\pi_ix}\otimes
{\pi_jf}\otimes {\pi_kg}\otimes \ldots )=0\ \ . $$
Hence $\mu$ is the required homomorphism, which is also onto.

 Conversely, by using the universal property of tensor product, we can define
$\lambda $ to be the homomorphism from $T(H_1,H_2,\ldots ,H_t)_{c+1}$ into
\\ $\ga_{c+1}(F)/\rho_{c+1}(S)\ga_{c+2}(F) $ given by
$$ \sum \underbrace{(h\otimes k\otimes \ldots )}_{(c+1)\
times}\longmapsto \underbrace{\prod [f,g,\ldots
]}_{wt.c+1}\rho_{c+1}(S)\ga_{c+2}(F)\ \ ,$$
 $$ with\ any\ bracketing\ \ \ \ \ \ \ \ bracketed\ the\ same\ way $$
where for $h\in H_i\ ,\ k\in H_j\ ,\ \ldots $, we pick
$f\in F_i,g\in F_j,\ldots ,$ such that $\pi_if=h\ ,\
\pi_jg=k\ ,\ldots $ . Clearly, this is a well-defined map since the
commutators on the right hand side are multilinear. One can easily see that
$\lambda$ is an epimorphism.

 Now the result follows , since $\mu \lambda $ and $\lambda \mu$ are the
identity maps on\\  $T(H_1,H_2,\ldots ,H_t)_{c+1}$ and
$\ga_{c+1}(F)/\rho_{c+1}(S)\ga_{c+2}(F)$ , respectively.\ \ $\Box$

In some aspect , the following theorem is a generalization of I.Schur [7,4],\\
J.Wiegold [8] , and M.R.R.Moghaddam [6] , where its proof follows from
Lemma 2.1 and Theorem 2.2 .\\
{\bf Theorem 2.3}

 Let $\prod_{i=1}^{t}\!^{\times}H_i$ be the direct product of finite cyclic
groups . Then by the above notation , the higher Schur-multiplicator of $G$
is as follows :
$$ {\cal N}_cM(\prod_{i=1}^{t}\!^{\times }H_i)\cong T(H_1,H_2,\ldots
,H_t)_{c+1}\ \ .$$

 Now we are ready to give an explicit formula for the higher Schur-multiplicator
of a finite abelian group with respect to the variety of nilpotent groups of
class at most $c\geq 1$ , ${\cal N}_c$ .

 Let $G$ be an arbitrary finite abelian group , then by the fundamental theorem
of finitely generated abelian groups , $G\cong {\bf Z}_{n_1}\oplus {\bf Z}_{n_2}
\oplus \ldots \oplus {\bf Z}_{n_k}$,\\ where $n_{i+1}|n_i$ for all $1\leq i\leq
k-1$ and $k\geq 2$ . ${\bf Z}_n^{(m)}$ will denote the direct product of $m$
copies of the cyclic group ${\bf Z}_n$ . Then with the above assumption, we
obtain the following theorem , which is a vast generalization of I.Schur (see
[4 or 7]).\\
{\bf Theorem 2.4}

 Let $G={\bf Z}_{n_1}\oplus {\bf Z}_{n_2}\oplus \ldots \oplus {\bf Z}_{n_k}$
be a finite abelian group . Then , for all $c\geq 1$ , the higher
Schur-multiplicator of $G$ is
$$ {\cal N}_cM(G)\cong {\bf Z}_{n_2}^{(b_2)}\oplus {\bf
Z}_{n_3}^{(b_3-b_2)}\oplus \ldots \oplus {\bf
Z}_{n_k}^{(b_k-b_{k-1})} \ \ \ \ ,$$
where $b_i$ is the number of basic commutators of weight $c+1$ on
$i$ letters.\\
{\bf Proof.}

Clearly ${\bf Z}_m\otimes {\bf Z}_n\cong {\bf Z}_{(m,n)}$ for all $m,n\in
{\bf N}$ , where $(m,n)$ is the greatest common divisor of $m$ and $n$ . Hence
$$ {\bf Z}_{m_1}\otimes {\bf Z}_{m_2}\otimes \ldots \otimes {\bf Z}_{m_k}\cong
{\bf Z}_{m_k}  \ , $$
when $m_{i+1}|m_i$ for all $1\leq i\leq k-1$ . Thus
$$ T({\bf Z}_{n_1},{\bf Z}_{n_2})_{c+1}\cong {\bf Z}_{n_2} \ \ .$$

 Now, by induction hypothesis assume
$$ T({\bf Z}_{n_1},{\bf Z}_{n_2},\ldots ,{\bf Z}_{n_{k-1}})_{c+1}\cong
{\bf Z}_{n_2}^{(b_2)}\oplus {\bf Z}_{n_3}^{(b_3-b_2)}\oplus
\ldots \oplus {\bf Z}_{n_{k-1}}^{(b_{k-1}-b_{k-2})} \ \ . $$
Then we have
$$ T({\bf Z}_{n_1},{\bf Z}_{n_2},\ldots ,{\bf Z}_{n_k})_{c+1} =
   T({\bf Z}_{n_1},{\bf Z}_{n_2},\ldots ,{\bf Z}_{n_{k-1}})_{c+1}\oplus L\ , $$
where $L$ is the summation of all those tensor products of ${\bf Z}_{n_1},
{\bf Z}_{n_2},\ldots ,{\bf Z}_{n_k} $ corresponding to the basic commutators of
weight $c+1$ on $k$ letters which involve ${\bf Z}_{n_k}$ . Since $n_k|n_i$ for
all $1\leq i\leq k-1$ , all those tensor products are isomorphic to
${\bf Z}_{n_k}$ . So $L$ is the direct product of $(b_k-b_{k-1})$ copies
of ${\bf Z}_{n_k}$ .\\
Hence the result follows, by induction.\ \  $\Box$

\end{document}